\documentclass[10pt,oneside]{article}
\input amssymb.sty
\usepackage[cp1250]{inputenc}
\usepackage[T1]{fontenc}

\begin{document}
\newcommand{\qed}{\rule{1.5mm}{1.5mm}}
\newcommand{\proof}{\textit{Proof. }}
\newcommand{\ccon}{\rightarrowtail}
\newtheorem{theorem}{Theorem}[section]
\newtheorem{lemma}[theorem]{Lemma}
\newtheorem{remark}[theorem]{Remark}
\newtheorem{example}[theorem]{Example}
\newtheorem{corollary}[theorem]{Corollary}
\newtheorem{proposition}[theorem]{Proposition}
\newtheorem{claim}[theorem]{Claim}
\begin{center}
\vspace*{0mm}
{\LARGE\textbf{Approximation of analytic sets with proper projection by algebraic sets}\vspace*{3mm}}\\
{\large Marcin Bilski\footnote{M. Bilski: Institute of
Mathematics, Jagiellonian University, \L ojasiewicza 6, 30-348
Krak\'ow,
Poland. e-mail: Marcin.Bilski@im.uj.edu.pl }}\vspace*{10mm}\\
\end{center}
\begin{abstract}Let $X$ be an analytic subset of $U\times\mathbf{C}^n$
of pure dimension $k$ such that the projection of $X$ onto $U$ is
a proper mapping, where $U\subset\mathbf{C}^k$ is a Runge domain.
We show that $X$ can be approximated by algebraic sets.
\vspace*{2mm}\vspace*{0mm}
\end{abstract}
\section{Introduction}
\label{intro} The problem of approximation of analytic objects by
simpler algebraic ones appears in a natural way in complex
geometry and has attracted the attention of several mathematicians
(see \cite{Ar}, \cite{BoK}, \cite{BMT1}, \cite{DLS}, \cite{Fo},
\cite{Ku}, \cite{Lem}, \cite{TT1}, \cite{TT2}, \cite{TT3}).

The aim of the present paper is to show that every purely
$k$-dimensional\linebreak analytic subset of $U\times\mathbf{C}^n$
whose projection onto $U$ is a proper mapping, where $U$ is a
Runge domain in $\mathbf{C}^k,$ can be approximated by purely
$k$-dimensional\linebreak algebraic sets (see Theorem \ref{main}).
Here by a Runge domain we mean a domain of holomorphy
$U\subset\mathbf{C}^k$ such that every function
$f\in\mathcal{O}(U)$ can be uniformly approximated on every
compact subset of $U$ by polynomials in $k$ complex variables (cf.
\cite{Ho}, pp 36, 52, and note that in general a Runge domain need
not be connected). The approximation is expressed in terms of the
convergence of holomorphic chains i.e. analytic sets are treated
as holomorphic chains with components of multiplicity one (see
Section~ \ref{holchai}).

The starting point for our considerations is the fact that purely
$k$-dimensional analytic subsets of $U\times\mathbf{C}^n$ with
proper projection onto a Runge domain $U\subset\mathbf{C}^k$ can
be approximated by complex Nash sets. (We discussed various
methods of such approximation in \cite{B1}, \cite{B2}, \cite{B3},
\cite{B4}, \cite{B5}, \cite{B6}.) This allows to reduce the proof
of Theorem \ref{main} to the case where the approximated object is
a complex Nash set. The problem in the reduced case is solved by
Proposition~\ref{ndimensional}.

Every point of a purely $k$-dimensional analytic subset $X$ of
some open $\Omega\subset\mathbf{C}^m$ has a neighborhood
$U\subset\Omega$ such that, after a linear change of the
coordinates, the projection of $X\cap U$ onto an open ball in
$\mathbf{C}^k\times\{0\}^{m-k}$ is a proper mapping. Therefore
Theorem \ref{main} immediately implies that every analytic set can
be locally approximated by algebraic ones in the sense of
holomorphic chains (see Corollary \ref{co}).

An important question, strongly motivated by the fact that
algebraic approximation is one of central techniques used in
numerical computations, is when the presented methods lead to
effective algorithms. This problem will be discussed in a
forthcoming paper. Here let us only mention that for Nash
approximations of analytic sets we described a constructive method
in \cite{B4}.

Note that the convergence of positive chains appearing Theorem
\ref{main} is equivalent to the convergence of currents of
integration over the considered sets (see \cite{Ch}, pp. 141,
206-207). The organization of this paper is as follows. In
Section~ \ref{prel} preliminary material is presented whereas
Section~ \ref{mainproof} contains the proofs of the main results.
\section{Preliminaries}
\label{prel}
\subsection{Nash sets}\label{prelnash}
Let $\Omega$ be an open subset of $\mathbf{C}^n$ and let $f$ be a
holomorphic function on $\Omega.$ We say that $f$ is a Nash
function at $x_0\in\Omega$ if there exist an open neighborhood $U$
of $x_0$ and a polynomial
$P:\mathbf{C}^n\times\mathbf{C}\rightarrow\mathbf{C},$ $P\neq 0,$
such that $P(x,f(x))=0$ for $x\in U.$ A holomorphic function
defined on $\Omega$ is said to be a Nash function if it is a Nash
function at every point of $\Omega.$ A holomorphic mapping defined
on $\Omega$ with values in $\mathbf{C}^N$ is said to be a Nash
mapping if each of its components is a Nash function.

A subset $Y$ of an open set $\Omega\subset\mathbf{C}^n$ is said to
be a Nash subset of $\Omega$ if and only if for every
$y_0\in\Omega$ there exists a neighborhood $U$ of $y_0$ in
$\Omega$ and there exist Nash functions $f_1,\ldots,f_s$ on $U$
such that $$Y\cap U=\{x\in U: f_1(x)=\ldots=f_s(x)=0\}.$$

The fact from \cite{Tw} stated below explains the relation between
Nash and algebraic sets.
\begin{theorem}\label{rednash}
Let $X$ be an irreducible Nash subset of an open set
$\Omega\subset\mathbf{C}^n.$ Then there exists an algebraic subset
$Y$ of $\mathbf{C}^n$ such that $X$ is an analytic irreducible
component of $Y\cap\Omega.$ Conversely, every analytic irreducible
component of $Y\cap\Omega$ is an irreducible Nash subset of
$\Omega.$
\end{theorem}
\subsection{Convergence of closed sets and holomorphic chains}\label{holchai}
Let $U$ be an open subset in $\mathbf{C}^m.$ By a holomorphic
chain in $U$ we mean the formal sum $A=\sum_{j\in J}\alpha_jC_j,$
where $\alpha_j\neq 0$ for $j\in J$ are integers and
$\{C_j\}_{j\in J}$ is a locally finite family of pairwise distinct
irreducible analytic subsets of $U$ (see \cite{Tw2}, cp. also
\cite{Ba}, \cite{Ch}). The set $\bigcup_{j\in J}C_j$ is called the
support of $A$ and is denoted by $|A|$ whereas the sets $C_j$ are
called the components of $A$ with multiplicities $\alpha_j.$ The
chain $A$ is called positive if $\alpha_j>0$ for all $j\in J.$ If
all the components of $A$ have the same dimension $n$ then $A$
will be called an $n-$chain.

Below we introduce the convergence of holomorphic chains in $U$.
To do this we first need the notion of the local uniform
convergence of closed sets. Let $Y,Y_{\nu}$ be closed subsets of
$U$ for $\nu\in\mathbf{N}.$ We say that $\{Y_{\nu}\}$ converges to
$Y$ locally uniformly if:\vspace*{2mm}\\
(1l) for every $a\in Y$ there exists a sequence $\{a_{\nu}\}$ such
that $a_{\nu}\in Y_{\nu}$
and\linebreak\hspace*{7mm}$a_{\nu}\rightarrow a$
in the standard topology of $\mathbf{C}^m,$\\
(2l)  for every compact subset $K$ of $U$ such that $K\cap
Y=\emptyset$ it holds $K\cap Y_{\nu}=\emptyset$\linebreak\hspace*{6.3mm}for almost all $\nu.$\vspace*{1mm}\\
Then we write $Y_{\nu}\rightarrow Y.$ For details concerning the
topology of local uniform convergence see \cite{Tw2}.

We say that a sequence $\{Z_{\nu}\}$ of positive $n$-chains
converges to a positive $n$-chain $Z$ if:\vspace*{2mm}\\
(1c) $|Z_{\nu}|\rightarrow |Z|,$\\
(2c) for each regular point $a$ of $|Z|$ and each submanifold $T$
of $U$ of dimension\linebreak\hspace*{7mm}$m-n$ transversal to
$|Z|$ at $a$ such that $\overline{T}$ is compact and
$|Z|\cap\overline{T}=\{a\},$\linebreak\hspace*{7mm}we have
$deg(Z_{\nu}\cdot T)=deg(Z\cdot T)$ for almost
all $\nu.$\vspace*{2mm}\\
Then we write $Z_{\nu}\ccon Z.$ (By $Z\cdot T$ we denote the
intersection product of $Z$ and $T$ (cf. \cite{Tw2}). Observe that
the chains $Z_{\nu}\cdot T$ and $Z\cdot T$ for sufficiently large
$\nu$ have finite supports and the degrees are well defined.
Recall that for a chain $A=\sum_{j=1}^d\alpha_j\{a_j\},$
$deg(A)=\sum_{j=1}^d\alpha_j$).

\subsection{Normalization of algebraic sets} Let us recall that
every affine algebraic set, regarded as an analytic space, has an
algebraic normalization (see \cite{Lo}, p. 471). Therefore (in
view of the basic properties of normal spaces, see \cite{Lo}, pp.
337, 343) the following theorem, which will be useful in the proof
of the main result, holds true.
\begin{theorem}\label{normal}Let $\tilde{Y}$ be an algebraic subset of
$\mathbf{C}^m.$ Then there are an integer $n$ and an algebraic
subset $Z$ of $\mathbf{C}^m\times\mathbf{C}^n$ with
$\pi(Z)=\tilde{Y},$ where
$\pi:\mathbf{C}^m\times\mathbf{C}^n\rightarrow\mathbf{C}^m$ is the
natural projection, satisfying the
following properties:\vspace*{2mm}\\
(01) $Z,$ regarded as an analytic space, is locally irreducible,\\
(02) $\pi|_{Z}:Z\rightarrow\mathbf{C}^m$ is a proper map,\\
(03) $\pi|_{Z\cap(\pi^{-1}(\tilde{Y}\setminus
Sing(\tilde{Y})))}:Z\cap(\pi^{-1}(\tilde{Y}\setminus
Sing(\tilde{Y})))\rightarrow\tilde{Y}$ is an injective map.
\end{theorem}
\subsection{Runge domains}\label{polyhed}
We say that $P$ is a polynomial polyhedron in $\mathbf{C}^n$ if
there exist polynomials in $n$ complex variables $q_1,\ldots,q_s$
and real constants $c_1,\ldots,c_s$ such that
$$P=\{x\in\mathbf{C}^n:|q_1(x)|\leq c_1,\ldots,|q_s(x)|\leq c_s\}.$$

The following lemma is a straightforward consequence of Theorem
2.7.3 and Lemma 2.7.4 form \cite{Ho}.
\begin{lemma}\label{hor1}Let $\Omega\subset\mathbf{C}^n$ be a Runge
domain. Then for every $\Omega_0\subset\subset\Omega$ there exists
a compact polynomial polyhedron $P\subset\Omega$ such that
$\Omega_0\subset\subset IntP.$
\end{lemma}
Theorem 2.7.3 from \cite{Ho} immediately implies the following
\begin{claim}\label{hor3}Let $P$ be any polynomial polyhedron in
$\mathbf{C}^n.$ Then $IntP$ is a Runge domain in $\mathbf{C}^n.$
\end{claim}

The following fact from \cite{Ho} (Theorem 2.7.7, p 55) will also
be useful to us.
\begin{theorem}\label{hor2}Let $f$ be a holomorphic function in a
neighborhood of a polynomially convex compact set
$K\subset\mathbf{C}^n.$ Then $f$ can be uniformly approximated on
$K$ by polynomials in $n$ complex variables.
\end{theorem}

\section{Approximation of analytic sets}\label{mainproof}
The following theorem is the main result of this paper.
\begin{theorem}\label{main}Let $U$ be a Runge domain in
$\mathbf{C}^k$ and let $X$ be an analytic subset of
$U\times\mathbf{C}^n$ of pure dimension $k$ with proper projection
onto $U.$ Then there is a sequence $\{X_{\nu}\}$ of algebraic
subsets of $\mathbf{C}^k\times\mathbf{C}^n$ of pure dimension $k$
such that $\{X_{\nu}\cap(U\times\mathbf{C}^n)\}$ converges to $X$
in the sense of holomorphic chains.
\end{theorem}

The proof of Theorem \ref{main} is based on two results. Firstly,
every purely dimensional analytic set with proper and surjective
projection onto a Runge domain can be approximated by Nash sets
(see Theorem \ref{nashap}). Secondly,  every complex Nash set with
proper projection onto a Runge domain can be approximated by
algebraic sets as stated in the following
\begin{proposition}\label{ndimensional}Let $Y$ be a Nash subset of
$\Omega\times\mathbf{C}$ of pure dimension $k<m,$ with proper
projection onto $\Omega,$ where $\Omega$ is a Runge domain in
$\mathbf{C}^{m-1}.$ Then there is a sequence $\{Y_{\nu}\}$ of
algebraic subsets of $\mathbf{C}^{m-1}\times\mathbf{C}$ of pure
dimension $k$ such that $\{Y_{\nu}\cap(\Omega\times\mathbf{C})\}$
converges to $Y$ in the sense of holomorphic chains.
\end{proposition}

\hspace*{-5.3mm}\textit{Proof of Proposition \ref{ndimensional}.}
Let $l$ be a positive integer and let $||\cdot||_{l}$ denote a
norm in $\mathbf{C}^l.$ Put
$B_l(r)=\{x\in\mathbf{C}^l:||x||_l<r\}.$ For any analytic subset
$X$ of an open subset of $\mathbf{C}^l$ let $X_{(q)}$ denote the
union of all $q$-dimensional irreducible components of $X.$

To prove the proposition it is clearly sufficient to show that for
every open $\Omega_0\subset\subset\Omega$ and
for every real number $r>0$ the following holds:\vspace*{2mm}\\
(*) there exists a sequence $\{Y_{\nu}\}$ of purely
$k$-dimensional algebraic subsets of\linebreak\hspace*{4.5mm}
$\mathbf{C}^{m-1}\times\mathbf{C}$ such that
$\{Y_{\nu}\cap(\Omega_0\times B_1(r))\}$ converges to
$Y\cap(\Omega_0\times B_1(r))$ in\linebreak\hspace*{4.9mm} the
sense of chains.\vspace*{2mm}

Fix an open relatively compact subset $\Omega_0$ of $\Omega$ and a
real number $r>0.$ Let
$\tilde{\pi}:\mathbf{C}^{m-1}\times\mathbf{C}\rightarrow\mathbf{C}^{m-1}$
denote the natural projection.
\begin{claim}\label{claim1} There exists a purely $k$-dimensional algebraic subset $\tilde{Y}$
of $\mathbf{C}^{m-1}\times\mathbf{C}$ such that
$Y\cap(\Omega_0\times B_1(r))$ is the union of some of the
analytic irreducible components of $\tilde{Y}\cap(\Omega_0\times
B_1(r)).$ Moreover, the mapping
$\tilde{\pi}|_{\tilde{Y}}:\tilde{Y}\rightarrow \mathbf{C}^{m-1}$
may be assumed to be proper.
\end{claim}

\hspace*{-5.3mm}\textit{Proof of Claim \ref{claim1}.} By Lemma
\ref{hor1} we can fix a compact polynomial polyhedron
$P\subset\Omega$ such that
$\Omega_0\subset\subset\Gamma\subset\subset\Omega,$ where
$\Gamma=IntP.$ The complex manifold
$Reg_{\mathbf{C}}(Y\cap(\Gamma\times\mathbf{C}))$ is a
semi-algebraic subset of $\mathbf{R}^{2m},$ hence it has a finite
number of connected components. Consequently,
$Y\cap(\Gamma\times\mathbf{C})$ has finitely many analytic
irreducible components. Therefore by Theorem~\ref{rednash} there
exists a purely $k$-dimensional algebraic subset $\tilde{Y}$ of
$\mathbf{C}^{m-1}\times\mathbf{C}$ such that
$Y\cap(\Gamma\times\mathbf{C})$ is the union of some of the
analytic irreducible components of
$\tilde{Y}\cap(\Gamma\times\mathbf{C}).$ Then, clearly,
$Y\cap(\Omega_0\times B_1(r))$ is the union of some of the
analytic irreducible components of $\tilde{Y}\cap(\Omega_0\times
B_1(r)).$

If the mapping
$\tilde{\pi}|_{\tilde{Y}}:\tilde{Y}\rightarrow\mathbf{C}^{m-1}$ is
proper then the proof is completed. Otherwise, using the facts
that
$Y\cap(\Gamma\times\mathbf{C})\subset\tilde{Y}\cap(\Gamma\times\mathbf{C})$
and $\Omega_0\subset\subset\Gamma,$ we show that there are a
$\mathbf{C}$-linear isomorphism
$\Phi:\mathbf{C}^m\rightarrow\mathbf{C}^m,$ a Runge domain
$\Omega_1$ in $\mathbf{C}^{m-1},$ and a real number $s>0$ such
that
the following hold:\vspace*{2mm}\\
(a) the projection of
$\Phi(\tilde{Y})\subset\mathbf{C}^{m-1}\times\mathbf{C}$ onto $\mathbf{C}^{m-1}$ is a proper mapping,\\
(b) $\Phi(\Omega_0\times
B_1(r))\subset\Omega_1\times B_1(s),$\\
(c) $\Phi(Y)\cap(\Omega_1\times B_1(s))$ is a Nash subset of
$\Omega_1\times\mathbf{C}$ whose projection onto
$\Omega_1$\linebreak\hspace*{5.8mm}is a proper
mapping,\\
(d) $\Phi(Y)\cap(\Omega_1\times B_1(s))$ is the union of some of
the irreducible components
of\linebreak\hspace*{5.8mm}$\Phi(\tilde{Y})\cap(\Omega_1\times
B_1(s)).$\vspace*{2mm}

If there exists a sequence $\{Z_{\nu}\}$ of purely $k$-dimensional
algebraic subsets of $\mathbf{C}^{m-1}\times\mathbf{C}$ such that
$\{Z_{\nu}\cap(\Omega_1\times B_1(s))\}$ converges to
$\Phi(Y)\cap(\Omega_1\times B_1(s))$ in the sense of chains, then
$Y\cap(\Omega_0\times B_1(r))$ is approximated, in view of (b), by
$\{\Phi^{-1}(Z_{\nu})\cap(\Omega_0\times B_1(r))\}.$ Moreover, (c)
implies that $\Omega_1$ and $\Phi(Y)\cap(\Omega_1\times B_1(s))$
taken in place of $\Omega$ and $Y$ respectively satisfy the
hypotheses of Proposition~\ref{ndimensional}. Since, in view of
(a) and (d), $\Phi(\tilde{Y})$ is a purely $k$-dimensional
algebraic subset of $\mathbf{C}^{m-1}\times\mathbf{C}$ with proper
projection onto $\mathbf{C}^{m-1},$ containing
$\Phi(Y)\cap(\Omega_1\times B_1(s)),$ the proof of the claim is
completed provided there are $\Phi,$ $\Omega_1$ and $s$ satisfying
(a), (b), (c) and (d).

Take $\Omega_1$ to be any Runge domain in $\mathbf{C}^{m-1}$ with
$\Omega_0\subset\subset\Omega_1\subset\subset\Gamma.$ (The
existence follows by Lemma \ref{hor1} and Claim \ref{hor3}.) Now,
since $dim(\tilde{Y})=k<m,$ by the Sadullaev theorem (see
\cite{Lo}, p. 389), the set $\mathcal{S}_{\tilde{Y}}$ of
one-dimensional linear subspaces $l$ of $\mathbf{C}^{m}$ such that
the projection of $\tilde{Y}$ along $l$ onto the orthogonal
complement $l^{\bot}$ of $l$ in $\mathbf{C}^m$ is proper, is open
and dense in the Grassmannian $\mathbf{G}_1(\mathbf{C}^m).$
Consequently, for every $\varepsilon>0$ there is
$l\in\mathcal{S}_{\tilde{Y}}$ so close to
$\{0\}^{m-1}\times\mathbf{C}$ that there is a $\mathbf{C}$-linear
isomorphism
$\Phi_{\varepsilon}:\mathbf{C}^m\rightarrow\mathbf{C}^m$
transforming $l,l^{\bot}$ onto $\{0\}^{m-1}\times\mathbf{C},$
$\mathbf{C}^{m-1}\times\{0\}$ respectively such that
$$||\Phi_{\varepsilon}-Id_{\mathbf{C}^m}||<\varepsilon,$$ where
$Id_{\mathbf{C}^m}$ is the identity on $\mathbf{C}^m.$

Clearly, $\Phi=\Phi_{\varepsilon}$ satisfies (a) (for every
$\varepsilon>0$). Now, by the facts that $Y$ is a Nash subset of
$\Omega\times\mathbf{C}$ such that
$\tilde{\pi}|_{Y}:Y\rightarrow\Omega$ is a proper map and
$\Omega_1\subset\subset\Omega,$ there is a real number $s>r$ such
that $$(\overline{\Omega_1}\times\partial B_1(s))\cap
Y=\emptyset.$$ This implies that $\Phi=\Phi_{\varepsilon}$
satisfies (c), if $\varepsilon$ is sufficiently small. Next, the
facts $\Omega_0\subset\subset\Omega_1$ and $s>r$ imply that
$\Phi=\Phi_{\varepsilon}$ satisfies (b), for small $\varepsilon.$
Finally, by $\Omega_1\subset\subset\Gamma$ we get
$\Phi_{\varepsilon}^{-1}(\Omega_1\times
B_1(s))\subset\Gamma\times\mathbf{C},$ if $\varepsilon$ is small
enough. Therefore $$\Phi_{\varepsilon}(Y)\cap(\Omega_1\times
B_1(s))\subset\Phi_{\varepsilon}(\tilde{Y})\cap(\Omega_1\times
B_1(s)),$$ which easily implies that (d) is satisfied with
$\Phi=\Phi_{\varepsilon}.$ Thus the
proof is completed.\qed\vspace*{2mm}\\
\textit{Proof of Proposition \ref{ndimensional} (continuation)}.
By Theorem \ref{normal} there are an integer $n$ and a locally
irreducible (regarded as an analytic space), purely
$k$-dimensional algebraic subset $Z$ of
$\mathbf{C}^m\times\mathbf{C}^n$ such that the restriction
$\pi|_{Z}$ of the natural projection
$\pi:\mathbf{C}^m\times\mathbf{C}^n\rightarrow\mathbf{C}^m$ is a
proper mapping, $\pi(Z)=\tilde{Y}$ and
$$\pi|_{Z\cap(\pi^{-1}(\tilde{Y}\setminus
Sing(\tilde{Y})))}:Z\cap(\pi^{-1}(\tilde{Y}\setminus
Sing(\tilde{Y})))\rightarrow\tilde{Y}$$ is an injective mapping.

We may assume that $(\tilde{Y}\setminus
Y)\cap(\Omega\times\mathbf{C})\neq\emptyset,$ because otherwise
$$Y\cap(\Omega_0\times B_1(r))=\tilde{Y}\cap(\Omega_0\times
B_1(r))$$ and one can take $Y_{\nu}=\tilde{Y},$ for every
$\nu\in\mathbf{N}.$

Now, the Nash subsets $E$ and $F$ of
$\Omega\times\mathbf{C}\times\mathbf{C}^n$ defined by
$$E=(Z\cap\pi^{-1}(Y))_{(k)} \mbox{ and }
F=(Z\cap\pi^{-1}(\overline{\tilde{Y}\setminus
Y})\cap(\Omega\times\mathbf{C}\times\mathbf{C}^n))_{(k)},$$ where
the closure is taken in $\Omega\times\mathbf{C},$ satisfy $E\cap
F=\emptyset.$ Indeed, if there exists some $a\in E\cap F$, then
$Z$ (regarded as an analytic space) is not irreducible at $a$
because
$$Z\cap(\Omega\times\mathbf{C}\times\mathbf{C}^n)=E\cup F,$$
and $dim(E\cap F)<k.$ Consequently, the sets
$$\tilde{E}=E\cap(P\times\mathbf{C}\times\mathbf{C}^n)\mbox{ and
}\tilde{F}=F\cap(P\times\mathbf{C}\times\mathbf{C}^n)$$ also
satisfy $\tilde{E}\cap\tilde{F}=\emptyset,$ where $P\subset\Omega$
is a fixed compact polynomial polyhedron such that
$\Omega_0\subset P.$ (The existence of $P$ follows by Lemma
\ref{hor1}.)

By Claim \ref{claim1}, we may assume that the mapping
$\tilde{\pi}|_{\tilde{Y}}$ is proper. Then the mapping
$\hat{\pi}|_Z:Z\rightarrow\mathbf{C}^{m-1}$ is proper as well,
where $\hat{\pi}=\tilde{\pi}\circ\pi.$ This implies that both
$\tilde{E}$ and $\tilde{F}$ are compact. Moreover, the mapping
$(\hat{\pi},p)|_{Z}:Z\rightarrow\mathbf{C}^{m}$ is proper for
every polynomial
$p:\mathbf{C}^{m}\times\mathbf{C}^n\rightarrow\mathbf{C}.$

Now the idea of the proof is to find a sequence $\{p_{\nu}\}$ of
polynomials defined on $\mathbf{C}^m\times\mathbf{C}^n$ with the
following properties:\vspace*{2mm}\\
(0) $\{p_{\nu}|_{\tilde{E}}\}$ converges uniformly to the mapping
$(x_1,\ldots,x_m,y_1,\ldots,y_n)\mapsto x_m,$\\
(1) $\inf_{b\in \tilde{F}}|p_{\nu}(b)|>r$ for almost all
$\nu.$\vspace*{2mm}\\
Then we show that the sequence $\{(\hat{\pi},p_{\nu})(Z)\}$ of
purely $k$-dimensional algebraic subsets of $\mathbf{C}^m$ is as
required in the condition (*):
$\{(\hat{\pi},p_{\nu})(Z)\cap(\Omega_0\times B_1(r))\}$ converges
to $Y\cap(\Omega_0\times B_1(r))$ in the sense of chains.

\begin{claim}\label{claim2}
There exists a sequence $\{p_{\nu}\}$ of polynomials in $m+n$
complex variables, satisfying (0) and (1).
\end{claim}
\textit{Proof of Claim \ref{claim2}.} Firstly, by the fact that
$\tilde{E}\cap\tilde{F}=\emptyset,$ there is an open subset $U$ of
$\mathbf{C}^m\times\mathbf{C}^n$ such that $U=U_1\cup U_2,$ where
$U_1, U_2$ are open subsets of $\mathbf{C}^m\times\mathbf{C}^n,$
$U_1\cap U_2=\emptyset$ and $\tilde{E}\subset U_1,
\tilde{F}\subset U_2.$

Secondly abbreviate $(x,y)=(x_1,\ldots,x_m,y_1,\ldots,y_n)$ and
note that the function $f:U\rightarrow\mathbf{C}$ defined by
$f(x,y)=x_m$ on $U_1$ and $f(x,y)=r+1$ on $U_2$ is holomorphic.

Thirdly observe that, since $Z$ is an algebraic subset of
$\mathbf{C}^m\times\mathbf{C}^n$ with proper projection onto
$\mathbf{C}^m,$ and $\pi(Z)=\tilde{Y}$ is an algebraic subset of
$\mathbf{C}^{m-1}\times\mathbf{C}$ with proper projection onto
$\mathbf{C}^{m-1},$ and $P$ is a compact polynomial polyhedron in
$\mathbf{C}^{m-1}$, the union $\tilde{E}\cup\tilde{F}$ is a
compact polynomial polyhedron (and hence a polynomially convex
compact set) in $\mathbf{C}^m\times\mathbf{C}^n.$ Indeed, there
are $s, \tilde{s}>0$ such that
$$\tilde{E}\cup\tilde{F}=Z\cap(P\times\overline{B_1(s)}\times
\overline{B_n(\tilde{s})})$$ and the righthand side of the latter
equation is clearly described by a finite number of inequalities
of the form $|Q(x,y)|\leq c,$ where $Q$ is a complex polynomial,
and $c$ is a real constant (possibly equal to zero).

Lastly, since $f$ is a holomorphic function in a neighborhood of a
polynomially convex compact set $\tilde{E}\cup\tilde{F},$ it is
sufficient to apply Theorem \ref{hor2} to obtain a sequence
$\{p_{\nu}\}$ of complex polynomials in $m+n$ variables converging
uniformly to $f$ on $\tilde{E}\cup\tilde{F}.$ Clearly, every such
sequence
satisfies (0) and (1).\qed\vspace*{2mm}\\
\textit{Proof of Proposition \ref{ndimensional} (end)}. Let
$\{p_{\nu}\}$ be a sequence of polynomials satisfying the
assertion of Claim \ref{claim2}. We check that the sequence
$\{Y_{\nu}\}$ defined by
$$Y_{\nu}=(\hat{\pi},p_{\nu})(Z),\mbox{ for }\nu\in\mathbf{N},$$
satisfies the condition (*), which is sufficient to complete the
proof of Proposition~\ref{ndimensional}.

Every $Y_{\nu}$ is a purely $k$-dimensional algebraic subset of
$\mathbf{C}^m$ because the mapping
$(\hat{\pi},p_{\nu}):\mathbf{C}^m\times\mathbf{C}^n\rightarrow\mathbf{C}^m$
is polynomial, its restriction $(\hat{\pi},p_{\nu})|_{Z}$ is
proper and $Z$ is a purely $k$-dimensional algebraic subset of
$\mathbf{C}^m\times\mathbf{C}^n.$ Hence it remains to check that
$\{Y_{\nu}\cap(\Omega_0\times B_1(r))\}$ converges to
$Y\cap(\Omega_0\times B_1(r))$ in the sense of chains.

To see that $\{Y_{\nu}\cap(\Omega_0\times B_1(r))\}$ converges to
$Y\cap(\Omega_0\times B_1(r))$ locally uniformly (c.f. (1l) and
(2l), Section \ref{holchai}) it is sufficient to observe that
$$Y_{\nu}\cap(\Omega_0\times B_1(r))=(\hat{\pi},p_{\nu})(\tilde{E})\cap(\Omega_0\times B_1(r)),$$ for
$\nu$ large enough, whereas
$$Y\cap(\Omega_0\times
B_1(r))=\pi(\tilde{E})\cap(\Omega_0\times B_1(r)).$$ The first
equality follows directly by (0), (1) and the definitions of
$\tilde{E}$ and $\tilde{F}.$ The latter one is an obvious
consequence of the definition of $\tilde{E}.$ Moreover, (0)
implies that $(\hat{\pi},p_{\nu})|_{\tilde{E}}$ converges to
$\pi|_{\tilde{E}}$ uniformly which in turn implies that both (1l)
and (2l) are satisfied.

To finish the proof, let us verify that
$\{Y_{\nu}\cap(\Omega_0\times B_1(r))\}$ and $Y\cap(\Omega_0\times
B_1(r))$ satisfy (2c) (of Section \ref{holchai}). Suppose that it
is not true. Then there exist a $k$-dimensional
$\mathbf{C}$-linear subspace $l$ of $\mathbf{C}^{m-1}\times\{0\}$
and open balls $C_1, C_2$ in $l, l^{\bot}$ respectively, where
$l^{\bot}$ denotes the orthogonal complement of $l$ in
$\mathbf{C}^m,$ such that $\overline{C_1+
C_2}\subset\Omega_0\times B_1(r)$ and the following
hold:\vspace*{2mm}\\
(a) $Y\cap(\overline{C_1}+\partial C_2)=\emptyset$\mbox{ and }$(\overline{\tilde{Y}\setminus Y})\cap\overline{(C_1+C_2)}=\emptyset,$\\
(b) every fiber of the projection of $Y\cap(C_1+C_2)$ onto $C_1$
is $1$-element,\\
(c) the generic fiber of the projection of $Y_{\nu}\cap(C_1+C_2)$
onto $C_1$ is at least\linebreak\hspace*{5.6mm}$2$-element for
infinitely many $\nu.$\vspace*{1.5mm}

The existence of $l\subset\mathbf{C}^m$ (and $C_1,$ $C_2$) as
above is a direct consequence of the assumption that (2c) does not
hold. Since the projection of
$\tilde{Y}\subset\mathbf{C}^{m-1}\times\mathbf{C}$ onto
$\mathbf{C}^{m-1}$ is a proper mapping, the subspace $l$ can be
chosen in such a way that it is contained in
$\mathbf{C}^{m-1}\times\{0\}.$

The conditions (a), (b) and the facts that
$\pi|_{Z\cap(\pi^{-1}(\tilde{Y}\setminus Sing(\tilde{Y})))}$ is
injective and $\pi|_{Z}$ is proper imply that
$$Z\cap(({C_1}+{C_2})\times\mathbf{C}^n)=graph(G),$$
where $G\in\mathcal{O}({C_1},C_2\times\mathbf{C}^n).$
Consequently, by (0) and the inclusion
$l\subset\mathbf{C}^{m-1}\times\{0\}$, for $\nu$ large enough,
$$T_{\nu}:=(\hat{\pi},p_{\nu})(graph(G))$$ is a $k$-dimensional
analytic subset of $C_1+C_2$ whose projection onto $C_1$ has
$1$-element fibers. Therefore, by (c), for infinitely many $\nu,$
there is an analytic subset $H_{\nu}$ of $C_1+C_2$ of pure
dimension $k$ such that $dim(T_{\nu}\cap H_{\nu})<k$ and
$T_{\nu}\cup H_{\nu}=Y_{\nu}\cap(C_1+C_2).$

On the other hand, for large $\nu,$ it follows that
$$H_{\nu}\subset(\hat{\pi},p_{\nu})(graph(G))=T_{\nu},$$ which clearly contradicts the fact that $dim(T_{\nu}\cap H_{\nu})<k.$
The inclusion holds because, as observed previously, for large
$\nu$ we have
$$H_{\nu}\subset Y_{\nu}\cap(C_1+ C_2)=
(\hat{\pi},p_{\nu})(\tilde{E})\cap(C_1+C_2).$$ Moreover, by (0),
(a) and the fact that $l\subset\mathbf{C}^{m-1}\times\{0\}$ it
holds
$$(\hat{\pi},p_{\nu})(\tilde{E})\cap(C_1+C_2)=(\hat{\pi},p_{\nu})(\tilde{E}\cap(({C_1}+{C_2})\times\mathbf{C}^n)),$$
which implies that $H_{\nu}\subset T_{\nu}$ as
$$\tilde{E}\cap(({C_1}+{C_2})\times\mathbf{C}^n)\subset
graph(G).$$

Thus $\{Y_{\nu}\cap(\Omega_0\times B_1(r))\}$ and
$Y\cap(\Omega_0\times B_1(r))$ satisfy (2c) and
the proof of Proposition~\ref{ndimensional} is completed.\qed\vspace*{2mm}\\
\textit{Proof of Theorem \ref{main}}. Let us first recall that for
analytic covers there exist Nash approximations (for details see
\cite{B2} or \cite{B4}, or \cite{B6}):
\begin{theorem}
\label{nashap}Let $U$ be a connected Runge domain in
$\mathbf{C}^k$ and let $X$ be an analytic subset of
$U\times\mathbf{C}^n$ of pure dimension $k$ with proper projection
onto $U.$ Then for every open relatively compact subset $V$ of $U$
there is a sequence of Nash subsets of $V\times\mathbf{C}^n$ of
pure dimension $k$ with proper projection onto $V$, converging to
$X\cap(V\times\mathbf{C}^n)$ in the sense of holomorphic chains.
\end{theorem}
Fix a Runge domain $U$ in $\mathbf{C}^k$ and an analytic subset
$X$ of $U\times\mathbf{C}^n$ of pure dimension $k$ with proper
projection onto $U.$ Clearly, in order to prove
Theorem~\ref{main}, it is sufficient to prove the following
\begin{claim}\label{claim3}For every open $V\subset\subset U$ there exists a sequence $\{X_{\nu}\}$ of
purely\linebreak $k$-dimensional algebraic subsets of
$\mathbf{C}^k\times\mathbf{C}^n$ such that
$\{X_{\nu}\cap(V\times\mathbf{C}^n)\}$ converges to
$X\cap(V\times\mathbf{C}^n)$ in the sense of chains.
\end{claim}
Let us prove the claim. Fix an open $V\subset\subset U.$ Since,
without loss of generality, $V$ can be replaced by a larger
relatively compact Runge subdomain of $U$ (cf.~Section
\ref{polyhed}) we may assume that $V$ is a Runge domain. By
Theorem \ref{nashap}, there is a sequence $\{Z_{\nu}\}$ of purely
$k$-dimensional Nash subsets of $V\times\mathbf{C}^n,$ with proper
projection onto $V,$ converging to $X\cap(V\times\mathbf{C}^n)$ in
the sense of chains.(Formally in Theorem~\ref{nashap}, $U$ is
assumed to be connected, but this assumption can be easily omitted
by treating every connected component of $U$ separately.)

For every $\nu$, by Proposition \ref{ndimensional} applied with
$Y=Z_{\nu},$ $\Omega=V\times\mathbf{C}^{n-1}$ and $m=n+k,$ there
is a sequence $\{Y_{\nu,\mu}\}$ of algebraic subsets of
$\mathbf{C}^k\times\mathbf{C}^n$ of pure dimension $k$ such that
$\{Y_{\nu,\mu}\cap(V\times\mathbf{C}^n)\}$ converges to $Z_{\nu}$
in the sense of chains. Clearly, there is a function
$\alpha:\mathbf{N}\rightarrow\mathbf{N}$ such that
$\{Y_{\nu,\alpha(\nu)}\cap(V\times\mathbf{C}^n)\}$ converges to
$X\cap(V\times\mathbf{C}^n).$ Thus the proofs of Claim
\ref{claim3} and Theorem \ref{main} are completed.\qed\\
\\
An immediate consequence of Theorem \ref{main} is the following
\begin{corollary}\label{co}Let $X$ be a purely $k$-dimensional analytic
subset of some open $\Omega\subset\mathbf{C}^m.$ Then for every
$a\in X$ there are an open neighborhood $U$ of $a$ in $\Omega$ and
a sequence $\{X_{\nu}\}$ of purely $k$-dimensional algebraic
subsets of $\mathbf{C}^m$ such that $\{X_{\nu}\cap U\}$ converges
to
$X\cap U$ in the sense of holomorphic chains.\\
\end{corollary}
\textbf{Acknowledgement}. I express my gratitude to Professor
Marek Jarnicki for information on polynomial approximation of
holomorphic functions. I also wish to thank Dr. Marcin Dumnicki
and Dr. Rafa\l$\mbox{ }$Czy\.{z} for helpful discussions.



\begin{thebibliography}{}
\bibitem{Ar}
Artin, M.: \emph{Algebraic approximation of structures over
complete local rings.} Publ. I.H.E.S. \textbf{36}, 23-58 (1969)
\bibitem{Ba}
Barlet, D.: \textit{Espace analytique r\'eduit des cycles
analytiques complexes compacts d'un espace analytique complexe de
dimension finie.} Fonctions de plusieurs variables complexes, II,
S\'em. Fran\c{c}ois Norguet, 1974-1975, Lecture Notes in Math.,
\textbf{482}, pp. 1-158, Springer, Berlin 1975
\bibitem{B1}
Bilski, M.: \textit{On approximation of analytic sets.}
Manuscripta Math. \textbf{114}, 45-60 (2004)
\bibitem{B2}
Bilski, M.: \textit{Approximation of analytic sets with proper
projection by Nash sets.} C.R. Acad. Sci. Paris, Ser. I
\textbf{341}, 747-750 (2005)
\bibitem{B3}
Bilski, M.: \textit{Approximation of analytic sets by Nash
tangents of higher order.} Math. Z. \textbf{256}, 705-716 (2007)
\bibitem{B4}
Bilski, M.: \textit{Algebraic approximation of analytic sets and
mappings.} J. Math. Pures Appl. \textbf{90}, 312-327 (2008)
\bibitem{B5}
Bilski, M.: \textit{Approximation of analytic sets along Nash
subvarieties.} Preprint
\bibitem{B6}
Bilski, M.: \textit{Approximation of sets defined by polynomials
with holomorphic coefficients.} Preprint
\bibitem{BoK}
Bochnak, J., Kucharz, W.: \emph{Approximation of holomorphic maps
by algebraic morphisms.} Ann. Polon. Math. \textbf{80}, 85-92
(2003)
\bibitem{BMT1}
Braun, R. W., Meise, R., Taylor, B. A.: \textit{Higher order
tangents to analytic varieties along curves.} Canad. J. Math.
\textbf{55}, 64-90 (2003)
\bibitem{Ch}
Chirka, E. M.: Complex analytic sets. Kluwer Academic Publ.,
Dordrecht-Boston-London 1989
\bibitem{DLS}
Demailly, J.-P., Lempert, L., Shiffman, B.: \textit{Algebraic
approximation of holomorphic maps from Stein domains to projective
manifolds.} Duke Math. J. \textbf{76}, 333-363 (1994)
\bibitem{Fo}
Forstneri\v c, F.: \textit{Holomorphic flexibility properties of
complex manifolds.} Amer. J. Math. \textbf{128}, 239-270 (2006)
\bibitem{Ho}
H\"{o}rmander, L.: An introduction to complex analysis in several
variables. North-Holland, Amsterdam-New York-Oxford-Tokyo 1991
\bibitem{Ku}
Kucharz, W.: \emph{The Runge approximation problem for holomorphic
maps into Grassmannians.} Math. Z. \textbf{218}, 343-348 (1995)
\bibitem{Lem}
Lempert, L.: \textit{Algebraic approximations in analytic
geometry.} Invent. Math. \textbf{121}, 335-354 (1995)
\bibitem{Lo}
\L ojasiewicz, S.: Introduction to complex analytic geometry.
Birkh\"{a}user Verlag Basel, 1991
\bibitem{TT1}
Tancredi, A., Tognoli, A.: \textit{Relative approximation theorems
of Stein manifolds by Nash manifolds.} Boll. Un. Mat. It.
\textbf{3-A,} 343-350 (1989)
\bibitem{TT2}
Tancredi, A., Tognoli, A.: \textit{On the extension of Nash
functions.} Math. Ann. \textbf{288}, 595-604 (1990)
\bibitem{TT3}
Tancredi, A., Tognoli, A.: \textit{On the relative Nash
approximation of analytic maps.} Rev. Mat. Complut. \textbf{11},
185-201 (1998)
\bibitem{Tw}
Tworzewski, P.: \textit{Intersections of analytic sets with linear
subspaces.} Ann. Sc. Norm. Super. Pisa \textbf{17,} 227-271 (1990)
\bibitem{Tw2}
Tworzewski, P.: \textit{Intersection theory in complex analytic
geometry.} Ann. Polon. Math., \textbf{62.2} 177-191 (1995)
\end{thebibliography}
\end{document}